\documentstyle[12pt]{article}
\newlength{\abstractwidth}
\renewcommand{\thefootnote}{\fnsymbol{footnote}}
\renewcommand{\thanks}[1]{\footnote{#1}} % Use this for footnotes
\newcommand{\starttext}{ \setcounter{footnote}{0}
\renewcommand{\thefootnote}{\arabic{footnote}}}

\newcommand{\be}{\begin{equation}}
\newcommand{\bea}{\begin{eqnarray}}
\newcommand{\eea}{\end{eqnarray}} \newcommand{\ee}{\end{equation}}
 
 \def\ba{\begin{eqnarray}}
\def\ea{\end{eqnarray}}

\def\tr{{\rm tr}}

\def\log{\,{\rm log}\,}
\def\exp{\,{\rm exp}\,}

\def\ge{\geq}
\def\le{\leq}

\def\p{\partial}

\def\[{{\bf [}}
\def\]{{\bf ]}}

\def\ddbar{i\p\bar\p}

\def\mathbb{\bf}

\def\eqref{\ref}

\begin{document}
\starttext \baselineskip=18pt \setcounter{footnote}{0}
\newtheorem{theorem}{Theorem}
\newtheorem{lemma}{Lemma}
\newtheorem{corollary}{Corollary}
\newtheorem{definition}{Definition}
\newtheorem{conjecture}{Conjecture}
\newtheorem{proposition}{Proposition}
\newcommand{\beqref}[1]{(\ref{#1})}

%% On $L^\infty$ estimates for complex Monge-Amp\`ere equation ON $L^\infty$ ESTIMATES FOR COMPLEX MONGE-AMP\`ERE EQUATIONS

\begin{center}
{\Large \bf On the modulus of  continuity of solutions to complex Monge-Amp\`ere equations
%%%%     On $L^\infty$ estimates for complex Monge-Amp\`ere equations
%On $L^\infty$ estimates for complex Monge-Amp\`ere equations
\footnote{Work supported in part by the National Science Foundation under grant DMS-1855947.}}

\medskip
\centerline{Bin Guo, Duong H. Phong, Freid Tong\footnote{F.T. is supported by Harvard's Center for Mathematical Sciences and Applications. }, and Chuwen Wang}

\medskip

\begin{abstract}

{\footnotesize In this paper, we prove a uniform and sharp estimate for the modulus of continuity of solutions to complex Monge-Amp\`ere equations, using the PDE-based approach developed by the first three authors in their approach to supremum estimates for fully non-linear equations in K\"ahler geometry. As an application,  we derive a uniform diameter bound for K\"ahler metrics satisfying certain Monge-Amp\`ere equations.}

\end{abstract}

\end{center}

\baselineskip=15pt
\setcounter{equation}{0}
\setcounter{footnote}{0}

\section{Introduction}
\setcounter{equation}{0}

The complex Monge-Amp\`ere equation has been extensively studied ever since Yau's seminal work on the solution of the Calabi conjecture \cite{Y}. Notably, assuming that the right hand side is in some Orlicz space, Kolodziej \cite{K} showed using pluripotential theory that the solution must be in $L^\infty$ (see also \cite{GPT} for a recent PDE-based proof), and in fact continuous. When the right hand side is in $L^q$ for some $q>1$, it is known that the solution must be H\"older continuous \cite{K1, DDK}. However, there are examples showing that H\"older continuity may not hold when the right hand side is not in $L^q$ for any $q>1$.
In general, when the solution is not H\"older continuous, its modulus of continuity is not known. For the complex Monge-Amp\`ere equation, the modulus of continuity of the solution is especially important, as it is closely related to essential geometric properties of the corresponding K\"ahler metric such as its diameter. In particular, uniform bounds for the diameter of the metric are needed for Gromov-Hausdorff convergence, and this requires in turn uniform bounds for the modulus of continuity. This is the problem which we shall solve in the present paper.

\medskip

Let $(X,\omega_0)$ be a compact K\"ahler manifold of complex dimension $n$. We consider the complex Monge-Amp\`ere equation
\begin{equation}\label{eqn:MA}
(\omega_0 + \ddbar u)^n = e^F \omega_0^n,\quad \omega: = \omega_0 + \ddbar u>0\mbox{ and  } \inf_X u = 1,
\end{equation}
where $F\in C^\infty(X,\mathbb{R})$ satisfies the compatibility condition $\int_X e^F \omega_0^n = \int_X \omega_0^n$. We shall use the following Orlicz norm for the right hand side $e^F$,
%The goal of this note is to derive a uniform continuity of the solution $u$ to \beqref{eqn:MA}, assuming the right hand side $e^F$ lies in some Orlicz space. We denote
$$\| e^F\|_{L^1(\log L)^p} = \int_X e^F |F|^p \omega_0^n.$$ 
%More precisely, we  shall show that
\begin{theorem}\label{thm:main}
Fix $p>n$. There exists a constant $C>0$ depending on $n, p, \omega_0, \| e^F\|_{L^1(\log L)^p}$ such that the following uniform estimate holds
\begin{equation}\label{eqn:main}
|u(x) - u(y)| \le \frac{C}{|\log d(x,y)|^\alpha}
\end{equation}
for any $x,y\in X$. Here $d(x,y)$ denotes the geodesic distance of the two points $x,y$ in the Riemannian manifold $(X,\omega_0)$, and $\alpha = \min \{\frac{p-n}{n}, \frac{p}{n+1}\}$.
\end{theorem}

The case of Riemann surfaces (i.e. $n=1$) shows that the solution $u$ to \beqref{eqn:MA} may fail to be (uniformly) H\"{o}lder continuous if $e^F\not\in L^q$ for any $q>1$. Theorem \ref{thm:main} says that the solution is still continuous with order $O(|\log d|^{-\alpha})$. This is a remarkable estimate in itself, as logarithmic moduli of continuity are rarely encountered, and this is the first example of it that we are aware of in the partial differential equations arising in K\"ahler geometry. The example at the end of Section \ref{section 3} implies the exponent $\alpha>0$ is sharp.  

\medskip

We remark that the estimate \beqref{eqn:main} continues to hold when the function $e^F$ is not smooth. This can be seen by a smoothing argument combined with the stability estimate of complex Monge-Amp\`ere equations \cite{K2, GPTa} and Theorem \ref{thm:main}. The continuity of $u$ when $e^F$ is not smooth had been obtained by Kolodziej \cite{K} using pluripotential theory and an argument by contradiction. Theorem \ref{thm:main} sharpens the continuity estimate of $u$ and provides a uniform control. 

\medskip

In the case $e^F\in L^q(X,\omega_0^n)$ for some $q>1$, it is known from \cite{K1, DDK} that the solution $u$ is H\"older continuous, i.e. $|u(x) - u(y)|\le C d(x,y)^a$ for some $a\in (0,1)$. Our proof of Theorem \ref{thm:main} can be easily modified to give a new and PDE-based proof of the H\"older continuity of $u$, in the spirit of the proof of $L^\infty$ estimates developed in \cite{GPT}. We can readily show in this manner that, if $e^F\in L^q$ for some $q>1$ and $q^* = \frac{q}{q-1}$ is the conjugate exponent of $q$, then $|u(x) - u(y)| \le C d(x,y)^{\alpha_0}$ for $\alpha_0 = \frac{2}{1+ (n+1) q^*}$, for some constant $C>0$ depending only on $n, q, \omega_0$ and $\| e^F\|_{L^q}$.  For the sake of interested readers, we provide a sketch of the proof in Section \ref{Holder}.

\smallskip

The complex Monge-Amp\`ere equation \beqref{eqn:MA} plays an important role in finding canonical K\"ahler metrics in complex geometry. It is natural to study the geometry of the K\"ahler metric $\omega_u = \omega_0+ \ddbar u$ satisfying \beqref{eqn:MA}, for example, its Ricci curvature, volume growth of geodesic balls, and diameter bound. Under some general assumptions on $\ddbar F$ and $e^F$, it is known by the work of Fu-Guo-Song \cite{FGS} that these geometric quantities are indeed bounded in some sense.  Without the assumption on $\ddbar F$, Y. Li \cite{Li} proves a diameter bound of $(X, \omega_u)$ if the function $e^F\in L^q(X,\omega_0^n)$ for some $q>1$. His proof requires the H\"older continuity of $u$ proved in \cite{K1,DDK} and Morrey's lemma.  With the uniform modulus of continuity estimate in Theorem \ref{thm:main}, we can generalize Li's result with a weaker assumption on $e^F$.

\begin{theorem}\label{thm:main 2}
Let $\omega_u$ be the solution to \beqref{eqn:MA}. Suppose $p>3n$, in other words, $\alpha>2$, then there exists a constant $C>0$ depending only on $n, p, \omega_0, \| e^F\|_{L^1(\log L)^p}$ such that
$${\mathrm{diam}}(X, \omega_u)\le C.$$
\end{theorem}
The diameter bound of K\"ahler metrics satisfying certain complex Monge-Amp\`ere equations or K\"ahler-Einstein type equations is a necessary ingredient in the study of degenerations of these metric spaces in the Gromov-Hausdorff sense \cite{S}. Theorem \ref{thm:main 2} provides a uniform diameter bound under a mild assumption on $e^F$, and we expect it to be useful in studying the geometry of complex Monge-Amp\`ere equations.

%Considering the progress in the study of degeneration of a family of complex Monge-Amp\`ere equations \cite{S, GPTa, GPTW}, it is natural to ask whether the modulus of continuity estimate in Theorem \ref{thm:main} still holds when the fixed K\"ahler metric $\omega_0$ is replaced by a family of degenerating K\"ahler metrics $\{\omega_t\}$.  The main difficulty lies in the construction of a ``uniform'' $\rho_\delta u$  as in \beqref{eqn:Demailly}. We shall return to this problem in a future work.

\medskip

\noindent{\em Convention:} we say a constant $C>0$ is {\em uniform} if it depends only on $n, p>n, \omega_0$ and $\| e^F\|_{L^1(\log L)^p}$. The constants $C$'s in different lines may not be the same, but are all uniform unless otherwise stated.

\section{Preliminaries}
\setcounter{equation}{0}

We collect some necessary background materials from \cite{DDK}. Let $\rho: \mathbb R_+\to \mathbb R_+$ be a smoothing kernel which is supported in $[0,1]$ and normalized to satisfy $\int_{\mathbb R_+} \rho(t) dt = 1$ and $\rho(t) = \mathrm{const}$ for $t\in [0,3/4]$. Given a function $u\in L^1(X)$ and $\delta\in (0,1)$, we define its $\delta$-regularization to be
\begin{equation}\label{eqn:Demailly}\rho_\delta u(z)  = \frac{1}{\delta^{2n}}\int_{\zeta\in T_z X} u(\exp_z(\zeta)) \rho(\delta^{-1} |\zeta|_{\omega_0}^2  ) dV_{\omega_0}(\zeta), \end{equation}
where $\exp_z: T_zX\to X$ is the exponential map of the Riemannian manifold $(X,\omega_0)$. $\rho_\delta u$ is a suitable weighted average of $u$ over the geodesic ball $B_{\omega_0}(z,\delta)$, so if $0\le u\in PSH(X,\omega_0)$, by mean value inequality $\rho_\delta u(z)$ control the maximum of $u$ over $B_{\omega_0}(z,\delta/2)$.
The following lemma is proved in \cite{DDK, BD}.
\begin{lemma}\label{lemma 1}
Let $u$ be an $L^1$ function in $PSH(X,\omega_0)$. Then
\begin{enumerate}
\item (\cite{BD}) There exists a constant $K>0$ depending on the curvature of $(X,\omega_0)$ such that
$t\mapsto \rho_t u(z) + Kt^2$ is monotone increasing for any $z\in X$.

\item (\cite{DDK}) There exists a constant $C>0$ depending on only $n, \omega_0$ such that
\begin{equation}\label{eqn:close}
\int_X |\rho_\delta u - u| \omega_0^n \le C \delta^2,\quad \forall \delta\in (0,1].
\end{equation}
\end{enumerate}
\end{lemma}
For a given small $c>0$, we define the Kiselman-Legendre transformation of $u$ as
\begin{equation}\label{eqn:Kiselman}u_{c,\delta}(z) = \inf_{t\in (0,\delta]} \{ \rho_t u(z) + Kt^2 - c \log \frac{t}{\delta} - K\delta^2     \}   \end{equation}
where $K>0$ is the constant in Lemma \ref{lemma 1}. By applying Kiselman's minimum principle it can be shown that (see \cite{D, BD}) for $u\in PSH(X,\omega_0)$
\begin{equation}\label{eqn:lower bound}\omega_0 + \ddbar u_{c,\delta} \ge - (A c + K\delta^2)\omega_0\end{equation}
where $-A$ is a lower bound of the bisectional curvature of the fixed K\"ahler metric $\omega_0$.

\section{Proof of Theorem \ref{thm:main}}\label{section 3}
\setcounter{equation}{0}
Let $u\in PSH(X,\omega_0)$ be the solution to the equation \beqref{eqn:MA}, where we normalize $u$ so that $u\ge 1$. We assume $p>n$ in this section. First recall the $L^\infty$ estimate of $u$ in \cite{K, GPT}.

\begin{lemma}\label{lemma L}
There exists a constant $C_0>0$ depending on $n, p, \|e^F\|_{L^1(\log L)^p},\omega_0$ such that
$$1\le u\le C_0,\quad \mbox{on }X.$$
\end{lemma}
We fix a small $\delta>0$. Let $c = \frac{1}{|\log \delta|^\alpha}$ for some $\alpha = \min( \frac{p-n}{n}, \frac{p}{n+1}  )>0$, and $U_\delta = u_{c,\delta}$ be the Kiselman-Legendre transformation of $u$ at the level $c$ as in \beqref{eqn:Kiselman}. From \beqref{eqn:lower bound} we get
$$\omega_0 + \ddbar U_\delta \ge - (Ac + K\delta^2) \omega_0 \ge - A' c \omega_0$$
for some $A' = A'(n,\omega_0)>0$. Hence the function $u_\delta = \frac{U_\delta}{1+ A' c}\le U_\delta$ belongs to $PSH(X,\omega_0)$, i.e. $\omega_0 + \ddbar u_\delta\ge 0$. We note from \beqref{eqn:Kiselman} and normalization of $u$ that $u_\delta$ is positive.

For any $s\ge 0$ we denote the set
\begin{equation}\label{eqn:Es}E_s = \{ u \le -2\delta + (1-r) u_\delta   - s \},\end{equation}
where $r = |\log \delta|^{-\frac{p}{n+1}}>0$ is a small constant.

\begin{lemma}\label{lemma E}
There is a constant $C_1>0$ depending on  $n, p, \|e^F\|_{L^1(\log L)^p},\omega_0$ such that
$$\int_{E_0} e^F \omega_0^n \le \frac{C_1}{|\log \delta|^p}.$$
\end{lemma}
\noindent{\em Proof.} We observe the following elementary inclusions of sets
$$E_0 = \{2\delta \le (1-r) u_\delta - u\}\subset \{2\delta \le  u_\delta -  u\} \subset  \{2\delta \le U_\delta - u\}  \subset  \{2\delta \le \rho_\delta u - u\} =:\Omega,$$
where the last inclusion follows from the fact that $U_\delta \le \rho_\delta u$. Thus it suffices to prove the lemma for the domain $\Omega$.

We define $v = \log \frac{\rho_\delta u - u}{\delta^{3/2}}$ as a function on $\Omega$. It is clear that $v \ge \log \frac{2}{\delta^{1/2}}>0$ on $\Omega$. Take a weight function $\eta(x) = (\log (1+x))^p$ on $\mathbb R_+$, and we apply the generalized Young's inequality with this weight. It follows that at any point $z\in \Omega$
\bea
v^p e^{F} \le \int_0^{e^F} \eta(x) dx + \int_0^{v^p} \eta^{-1}(y) dy \le e^F(1+ |F|)^p + v^p e^{v}. \nonumber
\eea
Integrating the above over $\Omega$, we get
$$\int_\Omega v^p e^F \omega_0^n\le C + \int_\Omega C |\log\delta|^p \frac{\rho_\delta u - u}{\delta^{3/2}} \omega^n_0\le C, $$
where in the last inequality we use (2) in Lemma \ref{lemma 1}. Since on $\Omega$, $v\ge \log 2 + \frac 1 2 |\delta| \ge \frac 12 |\log \delta|$, we conclude that
$$\frac 1 {2^p} |\log \delta|^p \int_\Omega e^F \omega_0^n \le \int_\Omega v^p e^F \omega_0^n\le C.$$
Then the lemma follows easily.

\medskip

We consider the auxiliary equations
\begin{equation}\label{eqn:aux}
(\omega_0+ \ddbar \psi_{s,k})^n = \frac{ \tau_k( -u + (1-r) u_\delta - 2 \delta - s  )   }{A_{s,k}} e^F \omega_0^n,\quad \sup_X \psi_{s,k} = 0,
\end{equation}
where $\tau_k(x): \mathbb R\to \mathbb R_+$ is a sequence of positive smooth functions converging decreasingly and pointwise to $x\cdot\chi_{\mathbb R_+}(x)$ on $\mathbb R$, and
\begin{equation}\label{eqn:Ask}
A_{s,k} = \int_{X}  \tau_k( -u + (1-r) u_\delta - 2 \delta - s  )  e^F \omega_0^n \to \int_{E_s}  ( -u + (1-r) u_\delta - 2 \delta - s  ) e^F \omega_0^n =:A_s.
\end{equation}
as $k\to\infty$. The equation \beqref{eqn:aux} admits a unique smooth solution by Yau's theorem \cite{Y}. As in \cite{GPT, GPTa, GPTW}, we aim to compare $\psi_{s,k}$ with the solution $u$ to \beqref{eqn:MA}. Consider the following test function
$$\Psi = -\varepsilon ( -\psi_{s,k} + \Lambda  )^{\frac{n}{n+1}} + [ -u + (1-r) u_\delta - 2 \delta - s   ],$$
where \begin{equation}\label{eqn:choice}
\varepsilon = ( \frac{n+1}{n}  )^{\frac{n}{n+1}} A_{s,k}^{\frac{1}{n+1}},\quad \Lambda = \frac{n}{n+1} \frac{A_{s,k}}{r^{n+1}}.
\end{equation}
We claim that $\Psi\le 0$ on $X$.

If the maximum point $x_{\max}$ of $\Psi$ lies in $X\backslash E_s^\circ$, by definition of $E_s$, it is clear that $\Psi(x_{\max})<0$. If $x_{\max} \in E_s^\circ$, then by maximum principle, at $x_{\max}$
\bea \nonumber
0 &\ge& \Delta_\omega \Psi  \ge  \frac{n\varepsilon}{n+1} (-\psi_{s,k} + \Lambda)^{-\frac{1}{n+1}} \Delta_\omega \psi_{s,k}  + (1-r) \Delta_\omega u_\delta - \Delta_\omega u\\
& \ge & \nonumber \frac{n\varepsilon}{n+1} (-\psi _{s,k}+ \Lambda)^{-\frac{1}{n+1}} \tr_\omega \omega_{\psi_{s,k}} + (1-r)\tr_\omega \omega_{u_\delta} - n + (r - \frac{n\varepsilon}{n+1} (-\psi_{s,k} + \Lambda)^{-\frac{1}{n+1}}) \tr_{\omega}\omega_0\\
& \ge & \nonumber \frac{n^2\varepsilon}{n+1} (-\psi_{s,k} + \Lambda)^{-\frac{1}{n+1}}\Big (\frac{\omega^n_{\psi_{s,k}}}{\omega^n} \Big)^{1/n} - n\\
& =  & \nonumber \frac{n^2\varepsilon}{n+1} (-\psi_{s,k} + \Lambda)^{-\frac{1}{n+1}}\Big ( \frac{ \tau_k( -u + (1-r) u_\delta - 2 \delta - s  )   }{A_{s,k}}  \Big)^{1/n} - n\\
& \ge   & \nonumber \frac{n^2\varepsilon}{n+1} (-\psi_{s,k} + \Lambda)^{-\frac{1}{n+1}}\Big ( \frac{ -u + (1-r) u_\delta - 2 \delta - s   }{A_{s,k}}  \Big)^{1/n} - n
\eea
where we denote $\omega_{\varphi} = \omega_0+ \ddbar \varphi$ for a function $\varphi\in PSH(X,\omega_0)$, in the third line we used the arithmetic-geometric inequality, the fact that $u_\delta\in PSH(X,\omega_0)$ and the choices of $\Lambda$ and $\varepsilon$ in \beqref{eqn:choice}. It follows easily that $\Psi(x_{\max})\le 0$ in this case. Hence we have $\Psi\le 0$. From the definition of $\Psi$ we obtain that on $E_s$
$$\frac{( -u + (1-r) u_\delta - 2 \delta - s )^{{(n+1)/}{n}}  }{A_{s,k}^{1/n}   }\le  C (n) (-\psi_{s,k} + \frac{A_{s,k}}{r^{n+1}}   ).  $$
By the H\"ormander estimate \cite{Ti, H} we can find a small $\beta_0 = \beta_0(n,\omega_0)>0$ such that
\bea \nonumber
\int_{E_s} \exp\Big(\beta_0 \frac{( -u + (1-r) u_\delta - 2 \delta - s )^{{(n+1)/}{n}}  }{A_{s,k}^{1/n}   }  \Big) \omega_0^n & \le & \int_X \exp( -  C(n)\beta_0 \psi_{s,k} + C(n)\beta_0 \frac{A_{s,k}}{r^{n+1}}  )\\
&\le & C \exp( C \frac{A_{s,k}}{r^{n+1}}  ) .\label{eqn:alpha}
\eea
Letting $k\to\infty$ in \beqref{eqn:alpha} we obtain a Trudinger-type estimate
\begin{equation}\label{eqn:Tru}
\int_{E_s} \exp\Big(\beta_0 \frac{( -u + (1-r) u_\delta - 2 \delta - s )^{{(n+1)/}{n}}  }{A_{s}^{1/n}   }  \Big) \omega_0^n \le C \exp( C \frac{A_s}{r^{n+1}}).
\end{equation}
\begin{lemma}\label{lemma 4}
We have a uniform constant  $C>0$ independent of $\delta$ and $s$ such that
$$\frac{A_s}{r^{n+1}}\le C.$$

\end{lemma}
\noindent {\em Proof.} By the uniform $L^\infty$ bound of $u$ and $u_\delta$, we have
$$A_s = \int_{E_s} (-u + (1-r) u_\delta - 2 \delta - s) e^F \omega_0^n\le C \int_{E_s} e^F \omega_0^n \le C \int_{E_0} e^F \omega_0^n \le \frac{C}{|\log \delta|^p}.$$
The lemma follows from this and the choice of $r = |\log \delta|^{-p/(n+1)}$.

\medskip

Combined with the lemma above, \beqref{eqn:Tru} implies
\begin{equation}\label{eqn:True}
\int_{E_s} \exp\Big(\beta_0 \frac{( -u + (1-r) u_\delta - 2 \delta - s )^{{(n+1)/}{n}}  }{A_{s}^{1/n}   }  \Big) \omega_0^n \le C,
\end{equation}
for some uniform constant $C>0$. Given \beqref{eqn:True}, we can apply the generalized Young's inequality as in \cite{GPT} to conclude that
$$\int_{E_s} [ -u + (1-r) u_\delta - 2 \delta - s ]^{{(n+1)p/}{n}} e^F \omega_0^n\le C A_{s}^{p/n},   $$
where $C>0$ is independent of $s$ and $\delta$. By H\"{o}lder inequality we have
\bea \nonumber
A_s & \le & \Big (\int_{E_s} [ -u + (1-r) u_\delta - 2 \delta - s ]^{{(n+1)p/}{n}} e^F \omega_0^n \Big)^{\frac{n}{p(n+1)}} \Big(\int_{E_s} e^F \omega_0^n \Big)^{1/q}
\\
&\le & C A_s^{\frac 1{n+1}} \Big(\int_{E_s} e^F \omega_0^n \Big)^{1/q},
\eea
where $q = \frac{p(n+1)}{p(n+1) - n}$ is the conjugate exponent of $p(n+1)/n$. So we get
$$A_s\le C \Big(\int_{E_s} e^F \omega_0^n \Big)^{(1+n)/qn} =C \Big(\int_{E_s} e^F \omega_0^n \Big)^{1 + a_0}  $$
where $a_0 =\frac{p-n}{pn}>0$. If we denote $\phi(s) = \int_{E_s} ( -u + (1-r) u_\delta - 2 \delta - s ) e^F \omega_0^n$, then it follows easily that
\begin{equation}\label{eqn:iteration}
s' \phi(s+s')\le C_3 \phi(s)^{1+a_0},\quad \forall s\ge 0,\, s'\ge 0,
\end{equation} for some uniform constant $C_3>0$ independent of $\delta\in (0,1/2]$.
Next we apply a De Giorgi-type iteration argument (c.f. \cite{K, GPT, GPTa}). Choose a $\delta_0>0$ small which depends only on $n, p, \omega_0, \| e^F\|_{L^1(\log L)^p}$  such that for all $\delta\in (0,\delta_0]$, it follows from Lemma \ref{lemma E}
$$C_3 \phi(0)^{a_0} = C_3 \Big( \int_{E_0} e^F \omega_0^n \Big)^{a_0}\le \frac{C_3C_1^{a_0} }{|\log \delta|^{p a_0}}< \frac 1 2.$$
And $\phi(0)^{a_0}\le C |\log \delta|^{-p a_0}$. With this choice of $\delta_0$, by a simple iteration argument (e.g. \cite{GPT}), we get the set $E_s = \emptyset$ for all $s> S_\infty$, where
$$S_\infty\le \frac{2 C_3}{1-2^{-a_0}} \phi(0)^{a_0} \le \frac{C}{|\log \delta|^{p a_0}}.  $$
We thus conclude that
$$u_\delta - u \le 2 \delta + r u_\delta + \frac{C}{|\log \delta|^{p a_0}} ,\quad \mbox{ on }X.$$
From the definition of $u_\delta$ and $U_\delta$, we get that on $X$
\begin{equation}\label{eqn:final 1}U_\delta - u \le 2 \delta + r u_\delta + A' c u +\frac{C}{|\log \delta|^{p a_0}}. \end{equation}
At each point $z\in X$, there exists a $t_z\in (0,\delta]$ realizing the infimum of $U_\delta = u_{c,\delta}$ in the definition \beqref{eqn:Kiselman}. From \beqref{eqn:final 1} it holds that
$$\rho_{t_z} u + Kt_z^2 - u - c \log \frac{t_z}{\delta} - K \delta^2 \le 2 \delta + r u_\delta + A' c u +\frac{C}{|\log \delta|^{p a_0}}.$$
Also (1) in Lemma \ref{lemma 1} shows that $\rho_{t_z} u + Kt_z^2 - u \ge 0$. So
$$\log \frac{t_z}{\delta} \ge\frac{ - K\delta^2- 2\delta}{c} - \frac{r}{c} u_\delta - A' u - \frac{C}{|\log \delta|^{p a_0} c},$$
by the choice of $c = \frac{1}{|\log \delta|^{\alpha}}$ with $\alpha = \min(p a_0, \frac{p}{n+1})$, it follows that there exists a uniform constant $C>0$ such that $\log \frac{t_z}{\delta}\ge -C$, thus $t_z\ge \theta \delta$ for some uniform $\theta\in (0,1)$. Again by (1) in Lemma \ref{lemma 1} we have at $z\in X$
\bea \nonumber
\rho_{\theta \delta} u + K\theta^2 \delta^2 - u& \le & \rho_{t_z} u + Kt_z^2 - u\\
& \le & \nonumber K\delta^2 + c \log \frac{t_z}{\delta} + 2 \delta + r u_\delta + A' c u +\frac{C}{|\log \delta|^{p a_0}}\le  \frac{C}{|\log \delta|^{\alpha}}.
\eea
This yields that for any $z\in X$ and $\delta\in (0,\delta_0]$, $\rho_{\theta \delta} u(z) - u(z) \le \frac{C}{|\log \delta|^{\alpha}}$; or equivalently for any $\delta\in (0, \theta\delta_0]$
\begin{equation}\label{eqn:final}\rho_{\delta} u(z) - u(z) \le \frac{C_4}{|\log \delta|^{\alpha}},\end{equation} for some uniform constant $C_4>0$.

We are now ready to finish the proof of Theorem \ref{thm:main}. For any $\delta>0$, we denote
$$\bar u_\delta(z) = \max_{x\in \overline{B(z, \delta)}} u(x), $$ where $B(z,\delta)$ denotes the geodesic ball with center $z$ and radius $\delta$ in the Riemannian manifold $(X,\omega_0)$. We claim that there exists a uniform constant $C>0$ such that $\bar u_\delta (z) - u(z) \le \frac{C}{|\log \delta|^\alpha}$ for any $z\in X$ and $\delta \in (0, \theta\delta_0/\beta_n]$ for some $\beta_n>0$ sufficiently large depending only on $n$ and $\alpha$. This would be sufficient to prove the theorem. We follow the arguments in \cite{GKZ} closely. In the paragraphs below we shall assume $\delta_0>0$ is chosen to be smaller than the injectivity radius of $(X,\omega_0)$, so the exponential maps considered are diffeomorphisms on relevant domains.

\medskip

We denote $\Omega(\delta) = \sup_{z\in X}  (\bar u_\delta(z) - u(z))$. Define $A>0$ to be \beqref{eqn:A} which depends only on $n, p, \omega_0$ and $C_4>0$, then we claim that $\Omega(\delta) \le \frac{A}{|\log \delta|^\alpha} $ for any $\delta\in (0,\theta\delta_0/\beta_n]$. Suppose not, then there exists some $0<\delta'< \theta\delta_0/\beta_n$ such that $\Omega(\delta')> \frac{A}{|\log \delta'|^\alpha}$. We define
\begin{equation}\label{eqn:delta}
\delta: = \inf\{ 0<t < \frac{\theta\delta_0}{\beta_n} | ~ \Omega(s) \le \frac{A}{|\log s|^\alpha} \mbox{ for all }s\in [t,  \frac{\theta\delta_0}{\beta_n}]   \}.
\end{equation}
The existence of $\delta'$ implies that $\delta\ge \delta' >0$. Since $u$ is continuous and $X$ is compact, there exists $z_0\in X$ such that $\Omega(\delta) = \bar u_\delta(z_0) - u(z_0) = u(w_0) - u(z_0)$ for some $w_0\in  \overline{B(z_0,\delta)}$. From the definition of $\delta$ in \beqref{eqn:delta}, it follows that
\begin{equation}\label{eqn:cons} \Omega(\delta) = \frac{A}{|\log \delta|^\alpha},\quad  \mbox{and ~}  \Omega(s) \le \frac{A}{|\log s|^\alpha} \mbox{ for all }s\in [\delta,  \frac{\theta\delta_0}{\beta_n}]. \end{equation}
We fix a constant $b>1$ but close to $1$, and observe that  $d(x, w_0)\ge b \delta$ for any $x\in B(z_0,3b\delta) \backslash B(w_0, b\delta)$, hence by \beqref{eqn:delta}
\begin{equation}\label{eqn:delta 1}u(w_0) - u(x) \le \Omega(d(x, w_0)) \le \frac{A}{|\log d(x, w_0)|^\alpha} = \frac{|\log \delta|^\alpha}{|\log (6b \delta)|^\alpha} \Omega(\delta)\le \hat C \Omega(\delta),\end{equation}
where $\hat C>0$ is an upper bound of $\frac{|\log \delta|^\alpha}{|\log (6b \delta)|^\alpha} $ for all $\delta\in (0, \theta\delta_0/\beta_n]$, which is uniform. By taking $\beta_n$ large enough (depending only on $n$ and $ p$), we can choose $\hat C$ arbitrarily close to $1$, since   $\lim_{\delta\to 0} \frac{|\log \delta|^\alpha}{|\log (6b \delta)|^\alpha} = 1$. Thus we can assume that $\hat C < 1 + \frac{1}{100^n}$, say. \beqref{eqn:delta 1} yields that
\begin{equation}\label{eqn:delta 2}
u(x) \ge u(w_0) - \hat C \Omega(\delta),\quad \forall~ x\in  B(z_0,3b\delta) \backslash B(w_0, b\delta).
\end{equation}By the definition of $\rho_{3b\delta} u$ in \beqref{eqn:Demailly}, we have
\bea\nonumber
\rho_{3b\delta} u(z_0)&  = & \frac{1}{(3b\delta)^{2n}}\int_{\zeta\in T_{z_0} X} u(\exp_{z_0} (\zeta)) \rho( \frac{|\zeta|_\omega^2}{(3b\delta)^2}  ) dV_\omega(\zeta)\\
&  \ge  & \frac{1}{(3b\delta)^{2n}}\int_{F} u(\exp_{z_0} (\zeta)) \rho( \frac{|\zeta|_\omega^2}{(3b\delta)^2}  ) dV_\omega(\zeta)
 \nonumber+ (1-\varepsilon_0)({u(w_0) - \hat C \Omega(\delta)}) ,
\eea
where $F\subset T_{z_0}X$ is the inverse image of $B(w_0,b\delta)$ under $\exp_{z_0}: T_{z_0}\to X$, and $$\varepsilon_0 = \frac{1}{(3b\delta)^{2n}} \int_{ F}  \rho( \frac{|\zeta|_\omega^2}{(3b\delta)^2}  ) dV_\omega(\zeta)\in [0,1].$$ Note that  we can choose $\varepsilon_0\ge \frac{1}{4^{2n}}$.
By Gauss' Lemma, $|\zeta|_\omega^2 \le (b+1)^2 \delta^2$ for any $\zeta\in F$. By the choice of the kernel function $\rho$, we have $ \rho( \frac{|\zeta|_\omega^2}{(3b\delta)^2}  )  = \mathrm{const}$ for such $\zeta$.  So
\bea\nonumber
 \frac{1}{(3b\delta)^{2n}}\int_{F} u(\exp_{z_0} (\zeta)) \rho( \frac{|\zeta|_\omega^2}{(3b\delta)^2}  ) dV_\omega(\zeta) = \varepsilon_0 \frac{1}{\mu(B(w_0, b\delta))} \int_{B(w_0, b\delta)} u(z) d\mu(z),
\eea
where $d\mu = (\exp_{z_0})_* dV_{\omega(z_0)}$ is the pushforward of the ``Euclidean measure'' in $T_{z_0}X$ to $X$ under the exponential map $\exp_{z_0}$. We observe that $B(w_0, b\delta)$ can be viewed as a domain in the normal coordinates chart at $z_0$, and under this coordinates system, the measure $\mu$ differs from the Euclidean one by $C\delta$ (for some uniform $C=C(\omega_0)$). Moreover, $u + \varphi_{z_0}$ is pluri-subharmonic for some local potential $\varphi_{z_0}$ of $\omega_0$ which satisfies $|\phi_{z_0}|\le C \delta$ (e.g. consider $\varphi_{z_0} - \varphi_{z_0}(w_0)$ if necessary). Then by the standard mean-value inequality for subharmonic functions in Euclidean space, we get
$$\varepsilon_0 \frac{1}{\mu(B(w_0, b\delta))} \int_{B(w_0, b\delta)} u(z) d\mu(z)\ge \varepsilon_0 u(w_0) - C_5 \delta,$$
where $C_5 >0$ is a constant depending only on $n,\omega_0$. Combining the above we get
\bea\rho_{3b\delta} u(z_0) & \ge & u(w_0) - C_5 \delta - (1-\varepsilon_0) \hat C \Omega(\delta) \nonumber \\
&= & u(z_0) - C_5 \delta + (1- (1-\varepsilon_0) \hat C) \Omega(\delta)\nonumber\\
& \ge & u(z_0) - C_5 \delta + \frac{1}{100^n} \Omega(\delta), \nonumber
\eea
where in the last inequality, we use the choices of $\varepsilon_0\ge 4^{-2n}$ and $\hat C\le 1+ \frac{1}{100^n}$. Combined with \beqref{eqn:final}, this yields that
\begin{equation}\label{eqn:contradiction}\frac{C_4}{|\log 3b\delta|^\alpha} + C_5 \delta \ge 10^{-2n} \Omega(\delta)= 10^{-2n}\frac{A}{|\log \delta|^\alpha}.\end{equation}
If at the beginning we choose $A>0$ to be
\begin{equation}\label{eqn:A}A = 1 + \Big|\log \frac{\theta \delta_0}{\beta_n}\Big|^\alpha \Omega(\frac{\theta \delta_0}{\beta_n}) + \sup_{\delta\in (0, \theta\delta_0/\beta_n]} 10^{2n} |\log \delta|^\alpha (\frac{C_4}{|\log 3b\delta|^\alpha} + C_5 \delta  )\end{equation}
The proof of Theorem \ref{thm:main} is complete.

\medskip

\noindent{\bf Example 1.} Let ${\mathbf{D}}\subset{ \mathbb C}\subset {\mathbb{CP}}^1$ be the disk with radius $1/2$. Consider the function
$\varphi(z) = (-\log |z|^2)^{-a}$ for some $a>0$ defined on ${\mathbf D}$ where $z$ is the standard coordinate on $\mathbb C$. Straightforward calculations show that on ${\mathbf{D}}\backslash \{0\}$
$$\ddbar \varphi = a (a+1) \frac{i dz\wedge d\bar z}{|z|^2 (-\log |z|^2)^{a+2}} = e^{\tilde F}.$$
It is easy to see that $e^{\tilde F}\in L^1 (\log L)^p (\mathbf{D})$ for any $p< a + 1$. The exponent $\alpha$ in Theorem \ref{thm:main} is $\alpha = p-1$ in this case. This example shows that the exponent $\alpha$ is sharp. Moreover, $\varphi$ is not H\"older continuous for any exponent.

Though $\varphi$ is singular in this example, we can consider a regularization of $\varphi$, for example, $\varphi_\epsilon(z) =  (-\log(\epsilon + |z|^2))^{-a}$ for $\epsilon\to 0^+$ to get a smooth example. We can also glue $\varphi$ or $\varphi_\epsilon$ to the whole space $\mathbb{CP}^1$ to get an example on a compact K\"ahler manifold.

\medskip

\noindent{\bf Example 2.} Let $(X,\omega_0)$ be a compact K\"ahler manifold and $L\to X$ be a holomorphic line bundle over $X$. Suppose $s\in H^0(X,{\mathcal{ O}}_X(L))$ is a nonzero holomorphic section of ${\mathcal O}_X(L)$. Consider the following complex Monge-Amp\`ere equations (for $\epsilon\in (0,1]$)
$$(\omega_0 + \ddbar u_\epsilon)^n = \frac{C_\epsilon}{(\epsilon + |s|_h^2) (-\log(\epsilon +  |s|_h^2))^a} \omega_0^n,$$
where $h$ is a Hermitian metric on $L$ such that $|s|^2_h <1$, $C_\epsilon>0$ is a normalizing constant so that the equation is solvable, and $a>0$ is a constant. Note that the function on the RHS belongs to $L^1(\log L)^p$ for any $p< a+1$. So for any {\em fixed} $a> n - 1$, Theorem \ref{thm:main} implies a uniform estimate on the modulus of continuity of $u_\epsilon$. Note that this estimate still holds for $\epsilon = 0$, if we interpret the complex Monge-Amp\`ere equation in the Bedford-Taylor sense.

\section{H\"older estimates}
\label{Holder}

We provide now a sketch of the proof of the H\"older continuity stated in the Introduction. Since the proof is analogous to that of Theorem \ref{thm:main}, we shall only point out the major differences. We use the same notations and definitions as before.

\smallskip 
 We assume $e^F\in L^q(X,\omega_0^n)$ for some $q>1$. Denote $\alpha_0 = \frac{2}{1+ (n+1)q^*}$ where $q^* = \frac{q}{q-1}$. The sets $E_s$ in \beqref{eqn:Es} is replaced by  $E_s= \{ u\le -2 \delta^{\alpha_0} + (1-r) u_\delta -s \}$ for $s\ge 0$, and here $r = \delta^{(2-\alpha_0)/(n+1)q^*}$. With these choices, Lemma \ref{lemma E} holds as $\int_{E_0} e^F \omega_0^n \le C \delta ^{(2-\alpha_0)/q^*}$ by H\"older inequality. Lemma \ref{lemma 4} continues to hold by the choice of $r$. We can proceed exactly as before to conclude \beqref{eqn:iteration} with any $a_0<1$, since in this case we can take $p$ as large as we like. The same iteration argument gives that $u_\delta - u \le C \delta^{\alpha_0}$. In the choice of $c$ in $U_{c,\delta}$ we can take $c = \delta^{\alpha_0}$. This will give $\rho_{\theta \delta} u - u\le C \delta^{\alpha_0}$ for some uniform $\theta\in (0,1)$ and any $\delta\in (0,\delta_0]$ for some uniform $\delta_0\in (0,1]$. Then it suffices to finish the proof of H\"older continuity of $u$ by invoking the estimates in \cite{GKZ} or the direct arguments in our proof of Theorem \ref{thm:main}.

\medskip

\section{Geometric applications}
\setcounter{equation}{0}

In this section, we apply a trick from \cite{Li} to see that the uniform continuity of the solution to \beqref{eqn:MA} leads to the diameter bound of the K\"ahler metric $\omega_u = \omega_0+ \ddbar u$, where $u$ satisfies \beqref{eqn:MA}.

Recall a function $f: \mathbb R_+ \to \mathbb R_+$ is called {\em Dini continuous}, if $ \int_0^1 \frac{f(r)}{r}dr <\infty$. As before,  we denote $\Omega(r) = \sup_{d(x,y)\le r} |u (x) - u(y)|$ to be the modulus of continuity of $u$, which is the oscillation of $u$ over geodesic balls of radius $r$.

\begin{lemma}\label{lemma Li}
On the K\"ahler manifold $(X,\omega_0)$, let $u$ be a smooth and strictly $\omega_0$-PSH function. If $\sqrt{\Omega}$ is Dini continuous, then the diameter of the K\"ahler manifold $(X,\omega_u)$ is bounded by a constant depending on $\omega_0$ and $\int_0^1 \frac{\sqrt{\Omega(r)}}{r } dr$.
\end{lemma}
\noindent{\em Proof.} Since $(X,\omega_0)$ is compact, we can take a finite open cover $\{U_a\}_{a=1}^N$, where each $U_a$ is a bounded domain in $\mathbb C^n$, and without loss of generality we assume each $U_a$ is biholomorphic to the Euclidean ball $B_{\mathbb C^n}(0, 2)$ and $\{\frac 1 2U_a\}$ also covers $X$. It is clear that $\omega_0|_{U_a}$ is equivalent to $\omega_{\mathbb C^n}|_{U_a}$. For notational convenience, in the proof of this lemma we write $B_r(z) = B_{\mathbb C^n}(z, r)$ and $\omega_E = \omega_{\mathbb C^n}$.

\def\Xint#1{\mathchoice
{\XXint\displaystyle\textstyle{#1}}%
{\XXint\textstyle\scriptstyle{#1}}%
{\XXint\scriptstyle\scriptscriptstyle{#1}}%
{\XXint\scriptscriptstyle\scriptscriptstyle{#1}}%
\!\int}
\def\XXint#1#2#3{{\setbox0=\hbox{$#1{#2#3}{\int}$ }
\vcenter{\hbox{$#2#3$ }}\kern-.6\wd0}}
\def\ddashint{\Xint=}
\def\dashint{\Xint-}

We consider the function $\rho(z) = d_{\omega_u} (z, 0)$, which is a Lipschitz function.  We fix a cut-off function $\chi: \mathbb R_+\to [0,1]$ such that $\chi(x) = 1$ for $x\in [0, 1]$ and vanishes on $[2,\infty)$. Following \cite{Li}, we look at the integral of $|\nabla \rho|_{\omega_0}^2$. For any fixed $r<1$ and any $p\in\frac 1 2 U_a\cong B_1(0)$, we have
\bea \nonumber
\int_{B_r(p)} |\nabla \rho|^2_{\omega_E} \omega_E^n & \le & \int_{B_r(p)} |\nabla \rho|^2_{\omega_u}   (\tr_{\omega_E} \omega_{u} )\omega_E^n = \int_{B_r(p)}   (n + \Delta_{\omega_E} u) \omega_E^n \\
& \le & C r^{2n} +  \int_{B_{2r}(p)} \Delta_{\omega_E}\chi(\frac{d_E(z, p)}{r}) \cdot (u(z) - u(p)) \omega_E^n \nonumber\\
& \le & C r^{2n} + C r^{2n-2} \Omega(2r)  \nonumber,
\eea
where in the second line we apply the integration by parts. By Poincare inequality it follows that
\bea\label{eqn:Morrey}
\dashint_{B_r(p)} \Big( \rho - \rho_{r,p} \Big)^2 \omega_E^n\le r^2 \dashint_{B_r(p)} |\nabla \rho|^2_{\omega_E}\omega_E^n  \le C r^2 + C \Omega(r),
\eea
where $\dashint_{B_r(p)} f $ denotes the average of $f$ over the ball $B_r(p)$, $\rho_{r,p} = \dashint_{B_r(p)} \rho\omega_E^n $, and in the last inequality we have applied $\Omega(2r)\le 2 \Omega(r)$ which follows from the triangle inequality.

\medskip

We now follow closely the proof of the classical Morrey's lemma in PDE theory. %Replacing $r$ by $r/2$ in \beqref{eqn:Morrey}, we get$$\dashint_{B_{r/2}(p)} \Big( \rho - \rho_{r/2 ,p} \Big)^2 \omega_E^n\le r^2 \le C (r/2)^2 + C \Omega(r/2).$$% 
By H\"{o}lder inequality and \beqref{eqn:Morrey}
\begin{equation}\label{eqn:Morrey 1}
|\rho_{r,p} -\rho_{r/2,p}   | \le \dashint_{B_{r/2}(p)} | u(z) - \rho_{r, p}  |   \omega_E^n \le C r + C \sqrt{\Omega(r)}.
\end{equation}
We apply \beqref{eqn:Morrey 1} with $r = 2^{-j}$ for $j = 1,2,3,\cdots$. Then
\begin{equation}\label{eqn:Morrey 2}
|\rho_{2^{-j}, p} - \rho_{2^{-j-1},p} |\le C 2^{-j}  + C \Omega(2^{-j})^{1/2}.
\end{equation}
Under the assumption that $\sqrt{\Omega(r)}$ is {\em Dini continuous}, we see that  the series $\hat \rho = \sum_{j=1}^\infty (\rho_{2^{-j}, p} - \rho_{2^{-j-1},p} )$ converges absolutely, and $|\hat \rho|$ is uniformly bounded, since $\sum_j 2^{-j}$ converges and $\sum_{j} \Omega(2^{-j})^{1/2} \le 2 \int_0^1 \frac{\sqrt{\Omega(t)}}{t}dt<\infty$. By Lebesgue differentiation theorem it is clear that \begin{equation}\label{eqn:desired}d_{\omega_u}(p, 0) = \rho(p) =\lim_{j\to \infty} \rho_{2^{-j}, p} = \rho_{1/2,p} -  \hat \rho. \end{equation} To get the desired bound on $d_{\omega_u}(p, 0)$ it suffices to estimate $\rho_{1/2,p}$. To this end, we observe that the inequalities above are uniform for any $p\in B_1(0)$. In particular we can apply \beqref{eqn:Morrey 1} and  \beqref{eqn:Morrey 2} with $r = 3\cdot 2^{-j}$ for $j=1,2,\cdots$ and $p=0$ to conclude that
$$d_{\omega_u}(0,0) = 0 = \rho_{3/2, 0} - O(1)   $$
where $O(1)$ denotes a uniformly bounded constant. This gives the bound on $\rho_{3/2, 0 } = \dashint_{B_{3/2}} \rho$.  Finally for any $p\in B_1(0)$, we have $B_{1/2} (p) \subset B_{3/2}(0)$ by triangle inequality, hence
$$\rho_{1/2, p} = \dashint_{B_{1/2}(p)} \rho \omega_E^n \le C \dashint_{B_{3/2}(0)} \rho \omega_E^n = C \rho_{3/2, 0} $$
is uniformly bounded, as desired. Combined with \beqref{eqn:desired}, this gives the expected bound on $d_{\omega_u}(p, 0)$ for any $p\in B_{1}(0)$. Since finitely many these balls cover $(X,\omega_0)$, we get the diameter bound of $(X,\omega_u)$. The proof of the lemma is complete.

\medskip

\noindent{\em Proof of Theorem \ref{thm:main 2}}. Let $u$ be the solution to \beqref{eqn:MA}. Suppose $p>3n$, then
$$\alpha = \min\{\frac{p - n}{n}, \frac{p}{1+n}\}> 1.$$
Theorem \ref{thm:main} implies that $|u(x) - u(y)|\le \frac{C}{|\log d(x,y)|^\alpha}$ for any $x,y\in X$. So we have for the modulus of continuity of $u$, $\Omega(r)\le \frac{C}{|\log r|^\alpha}$. It is now elementary to see that
$$\int_0^{1/2} \frac{\sqrt{\Omega(r)}}{r} dr \le \int_{\log 2} ^\infty \frac{C}{t^{\frac{\alpha}{2} }} dt <\infty. $$
We can now apply Lemma \ref{lemma Li} to conclude the uniform diameter bound of $(X,\omega_u)$.

\medskip

\noindent {\bf Example.} In {\bf Example 2} at the end of Section \ref{section 3}, if $a> 3n - 1$, Theorem \ref{thm:main 2} implies a uniform diameter bound of the K\"ahler metrics $\omega_\epsilon = \omega_0+ \ddbar u_\epsilon$, which is independent of $\epsilon\in (0,1]$.

\bigskip

\noindent{\bf Acknowledgement:} Bin Guo would like to thank Professor Jian Song for many helpful and stimulating discussions on complex Monge-Amp\`ere equations.

\bigskip

%Email address:  phong@math.columbia.edu, tong@math.columbia.edu

\noindent Department of Mathematics \& Computer Science, Rutgers University, Newark, NJ 07102

\noindent bguo@rutgers.edu

\medskip

\noindent Department of Mathematics, Columbia University, New York, NY 10027

\noindent phong@math.columbia.edu, wang.chuwen@columbia.edu

\medskip
\noindent Center for Mathematical Sciences and Applications, Harvard University, Cambridge, MA 02138

\noindent ftong@cmsa.fas.harvard.edu

\end{document}